\documentclass[12pt]{amsart}
\usepackage{amsmath, amssymb}
\usepackage[applemac]{inputenc}
\usepackage{diagrams}
\newcommand{\bdi}{\begin{diagram}}
\newcommand{\edi}{\end{diagram}}
\diagramstyle[scriptlabels,height=8mm,width=8mm]
\usepackage{enumerate}

\theoremstyle{definition}
\newtheorem{dfn}{Definition}[section]
\newtheorem{rem}[dfn]{Remark}
\newtheorem{sit}[dfn]{}
\newtheorem{exa}[dfn]{Example}
\newtheorem{prob}[dfn]{Problem}

\theoremstyle{plain}
\newtheorem{prop}[dfn]{Proposition}
\newtheorem{thm}[dfn]{Theorem}
\newtheorem{lem}[dfn]{Lemma}
\newtheorem{cor}[dfn]{Corollary}

\newcommand{\brem}{\begin{rem}}
\newcommand{\erem}{\end{rem}}

\newcommand{\bexa}{\begin{exa}}
\newcommand{\eexa}{\end{exa}}

\newcommand{\bdefi}{\begin{dfn}}
\newcommand{\edefi}{\end{dfn}}

\newcommand{\bcor}{\begin{cor}}
\newcommand{\ecor}{\end{cor}}

\newcommand{\blem}{\begin{lem}}
\newcommand{\elem}{\end{lem}}

\newcommand{\bprop}{\begin{prop}}
\newcommand{\eprop}{\end{prop}}

\newcommand{\bprob}{\begin{prob}}
\newcommand{\eprob}{\end{prob}}

\newcommand{\bthm}{\begin{thm}}
\newcommand{\ethm}{\end{thm}}

\newcommand{\bsit}{\begin{sit}}
\newcommand{\esit}{\end{sit}}

\newcommand{\be}{\begin{equation}}
\newcommand{\ee}{\end{equation}}

\newcommand{\bproof}{\begin{proof}}
\newcommand{\eproof}{\end{proof}}

\def\ba{\begin{array}}
\def\ea{\end{array}}
\def\bea{\begin{eqnarray}}
\def\eea{\end{eqnarray}}
\def\bals#1\eals{\begin{align*}#1\end{align*}}
\def\bal#1\eal{\begin{align}#1\end{align}}

\def\bdet#1\edet{\left|\begin{matrix}#1\end{matrix}\right|}
\def\bsdet#1\esdet{\left|\begin{smallmatrix}#1\end{smallmatrix}\right|}

\def\bmat{\begin{pmatrix}}
\def\emat{\end{pmatrix}}

\def\bnum{\begin{enumerate}}
\def\enum{\end{enumerate}}

\def\bsmat{\left(\begin{smallmatrix}}
\def\esmat{\end{smallmatrix}\right)}

\def\be{\begin{equation}}
\def\ee{\end{equation}}

\def\hto{\hookrightarrow}

\newcommand{\Spec}{\operatorname{Spec}}

\newcommand{\Ext}{{\operatorname{Ext}}}

\newcommand{\sdim}{{\operatorname{s.gl.dim}}}
\newcommand{\wid}{{\operatorname{wd}}}
\newcommand{\pd}{{\operatorname{pd}}}

\def\fa{{\mathfrak a}}

\def\fm{{\mathfrak m}}

\def\fp{{\mathfrak p}}

\def\cC{{\mathcal C}}

\def\cF{{\mathcal F}}
\def\cG{{\mathcal G}}
\def\cH{{\mathcal H}}

\def\cM{{\mathcal M}}

\def\cO{{\mathcal O}}
\def\cP{{\mathcal P}}
\def\cQ{{\mathcal Q}}

\newcommand{\F}{{\mathbb F}}

\renewcommand{\and}{\quad\mbox{and}\quad}

\begin{document}

\title{Strong global dimension of commutative rings and schemes}
\dedicatory{To the memory of Dieter Happel}
\author{Ragnar-Olaf Buchweitz}
\address{Dept.\ of Computer and Mathematical Sciences, University of
Toronto Scarborough, Toronto, ON M1C 1A4, Canada}
\email{ragnar@utsc.utoronto.ca}

\author{Hubert Flenner}
\address{Fakult\"at f\"ur Mathematik der Ruhr-Universit\"at,
Universit\"atsstr.\ 150, Geb.\ NA 2/72, 44780 Bochum, Germany}
\email{hubert.flenner@rub.de}

\thanks{The first author was partly supported by NSERC grant
3-642-114-80 and by the Humboldt foundation.}


\begin{abstract}
The strong global dimension of a ring is the supremum of the length of perfect complexes that are indecomposable in the derived category. In this note we characterize the noetherian commutative rings that have finite strong global dimension. We also give a similar characterization for arbitrary noetherian schemes. 
\end{abstract}

\subjclass[2000]{13D05, 13D09, 14F05, 18G20}
\keywords{strong global dimension, perfect complex}
\maketitle

\section*{Introduction}

Following \cite{Il} by a perfect complex over a ring $A$ we mean  a complex $P^\bullet$ of $A$-modules that is quasiisomorphic to a finite complex 
$$
P^\bullet:\quad 0\to P^a\to\ldots \to P^{b-1}\to P^b\to 0.
$$
of finite projective $A$-modules. Such a complex is called {\em indecomposable} if it is not quasiisomorphic to a direct sum of two non-trivial perfect complexes. 
We recall from \cite{HZ1} that the {\em strong global dimension} $\sdim(A)$ of an algebra is the supremum of lengths $l:=b-a$ of indecomposable perfect complexes $P^\bullet$ as above.  For instance, it is an elementary fact that a field has strong global dimension 0 and that a Dedekind domain has strong global dimension 1. 

The notion of strong global dimension was first proposed by C.M. Ringel and studied  in \cite{Sk}, \cite{HZ1, HZ2}. In \cite{HZ1} it was shown that a finite dimensional algebra over a field has finite strong global dimension if and only if it is piecewise hereditary.  

In this note we will study this notion in the context of commutative algebra and algebraic geometry. A ring of finite strong global dimension will be called {\em strongly regular}\footnote{Note that our strongly regular rings are different from those considered in relation to von Neumann regular rings
\cite[p.\ 28]{Go}}. Our main result is that a noetherian ring $A$ is strongly regular if and only if it is a finite product of rings $A=\prod_{i=1}^s A_i$ where each $A_i$ is a field or a Dedekind domain. In a second part we deduce a similar result for  noetherian schemes.   

\section{Strongly regular rings}

\bdefi\label{1.1}
A ring $A$ is called {\em strongly $r$-regular} if every perfect complex over $A$ is quasiisomorphic to a direct sum of complexes of length $\le r$. If $A$ is strongly $r$-regular for some $r$ then it will be called {\em strongly regular}.  
\edefi

\brem\label{1.2}
1. By a result of \cite{Sk}, a hereditary ring $A$ is strongly $1$-regular. In particular this holds for products $A=\prod_{i=1}^sA_i$, where each $A_i$ is a Dedekind domain or a field. 

\bproof
Let us provide a simple argument. 
Let $F^\bullet$ be a perfect complex over $A$. As every $A$-module can be represented as perfect complex of length 1, it suffices to show that $F^\bullet$ is formal, i.e.\ quasiisomorphic to its cohomology. We may assume that $F^i$ vanishes for $i> 0$ and consider  $H:=H^0(F^\bullet)$ as a complex concentrated in degree 0. With $G^\bullet$  the kernel of 
the natural map $F^\bullet\to H$, 
we get an exact sequence of complexes 
\be\label{eq1.1}
0\to G^\bullet\to F^\bullet\to H\to 0\,.
\ee
Let us consider the spectral sequence 
$$
E^{p,q}_1=\Ext_A^q(H, H^p(G^\bullet)))\Rightarrow \Ext_A^{p+q}(H, G^\bullet)\,.
$$
If  $q\ne 0,1$ or $p\ge 0$ then $\Ext_A^q(H, H^p(G^\bullet))=0$. Thus, the group
$\Ext^{1}_A(H, G^\bullet)$ vanishes, whence the sequence \eqref{eq1.1} splits, i.e.\ $F^\bullet$ is quasiisomorphic to $ G^\bullet\oplus H$. A simple induction  now yields that $F^\bullet$ is formal, as required. 
\eproof

2. Every field or finite product of fields is $0$-regular. More generally, if $A$ is a semi-simple ring, i.e.\ every $A$-module is projective, then it is 0-regular. 

3. If a ring $A$ is strongly $r$-regular and noetherian then its global dimension  is $\le r$. If $A$ is furthermore commutative then $A$ is a regular ring of dimension $\le r$. 

\bproof 
In lack of an explicit  reference of this standard fact we provide a short argument. By a result of Auslander, see \cite[Theorem 9.12]{Ro}, a ring has finite global dimension at most $r$ if and only if every finite (even cyclic) $A$-module is of finite projective dimension $\le r$.
Let $M$ be a finite $A$-module and let us show that $M$ has a resolution of length $\le r$ by finite projectives. Consider a  free resolution $F^\bullet$ of $M$ by finite $A$-modules and let us cut this resolution in degree $-r-1$ to obtain a perfect complex 
$$
\tilde F^\bullet:\quad  0\to F^{-r-1}\to F^{-r}  \to \cdots\to  F^0\to 0\,.
$$
By assumption $\tilde F^\bullet $ is quasiisomorphic to $\bigoplus_{i=1}^s F^\bullet_{(i)}$, where each $F^\bullet_{(i)}$ is a perfect complex of length $\le r$. Let $I:=\{i: H^0(F^\bullet_{(i)})\ne 0\}$, so that 
$M=H^0(\tilde F^\bullet) =\bigoplus_{i\in I} H^0(F^\bullet_{(i)})$. By construction the complex 
$$
\bigoplus_{i\in I}\tilde F_{(i)}^\bullet
$$
is exact in all degrees $\ne 0, -r-1$ since it is up to quasiisomorphism a direct summand of $\tilde F^\bullet$. As its length is $\le r$ it is also exact in degree $-r-1$. Hence it constitutes a resolution of $M$ of length $\le r$ as desired. 
\eproof

\erem 

Our main result is as follows.  

\bthm\label{1.3}
A commutative noetherian ring $A$ is strongly regular if and only if it is a product $\prod_{i=1}^s A_i$, where each $A_i$ is a Dedekind domain or a field. 
\ethm

{\em For the remaining part of this section let $A$ always denote a commutative noetherian ring.}
To prove Theorem \ref{1.3} we will first treat the case that $A$ is local. In this case the key observation is contained in the next proposition. To formulate it we need a few notations. 

\bsit\label{1.4}
Let $(A,\fm)$ be a  local noetherian ring and $(x,y)$ a regular sequence in $A$. Consider the matrix 
$$
M:=\bmat xy & -x^2\\ y^2 &-xy\emat\,.
$$
Multiplying (row) vectors from the right with $M$ yields  a periodic infinite complex 
\bdi
\F:\quad \cdots & \rTo^{\cdot M} & A^2 &\rTo^{\cdot M} &A^2 &\rTo^{\cdot M} & A^2 &\rTo^{\cdot M} \cdots
\edi
We extract from this a complex of length $n$ 
\be\label{eq1}
\bdi
\F^\bullet_n:&0&\rTo & A^2 &\rTo^{\cdot M} & A^2 &\rTo^{\cdot M} &\cdots&\rTo^{\cdot M}  & A^2 &\rTo^{\cdot M}  & A^2 &\rTo^{\cdot \bsmat x\\-y\esmat } & A &\rTo 0\,,
\edi
\ee
where $\F_n^0\cong A$,  $\F_n^{i}\cong A^2$ for $-n\le i<0$ and $\F_n^{i}=0$ otherwise.
\esit

We claim the following.

\bprop\label{1.5}
The complex  $\F^\bullet_n$ is indecomposable.
In particular, $\F^\bullet_n$ is not quasiisomorphic to a direct sum of  perfect complexes of length $<n$.  
\eprop

To prove this we need a simple lemma.
We recall that a complex $P^\bullet$ of free $A$-modules is called {\em minimal} if $\partial (P^{i-1})\subseteq \fm P^i$. We have the following  standard fact.

\blem\label{1.6}
(1) Every perfect complex $P^\bullet$ of free $A$-modules is quasiisomorphic to a minimal one. 

(2) If two minimal perfect complexes of free $A$-modules are quasiisomorphic, i.e.\ isomorphic in the derived category $D_b(A)$, then they are isomorphic. 
\elem 

\bproof
(1) is well known, at least for free resolutions of modules, see e.g.\ \cite{Wei}. 
In lack of an explicit reference we give a short argument in the general case. 
Given a finite complex of free $A$-modules $P^\bullet$ let us show inductively that we can minimize it by splitting off acyclic free complexes. Let $k$ be the number of indices $i$ such that $\partial(P^{i-1})\not\subseteq \fm P^i$. We proceed by decreasing induction on $k$. If $k=0$ then $P^\bullet$ is minimal. 

Suppose that for some index $i$ we have $\partial(P^{i-1})\not\subseteq \fm P^i$ so that  there is a submodule  $Q^i\subseteq \partial(P^{i-1})$ which is a direct summand of $P^i$. Choose $Q^{i}$ maximal with this property and take a free submodule $Q^{i-1}\subseteq P^{i-1}$ which is isomorphic to $Q^i$ under $\partial$.  Then the 2-term complex 
$$
Q^\bullet:\quad 0\to Q^{i-1}\to Q^i\to 0
$$
is a free acyclic subcomplex of $P^\bullet$. The quotient $\bar P^\bullet=
P^\bullet/Q^\bullet$ is a  complex of free $A$-modules that is quasiisomorphic to $P^\bullet$ and  satisfies $\partial \bar P^{i-1}\subseteq \fm \bar P^i$ as required. This yields the induction step and thus (1).

(2) Let $P^\bullet$ and $Q^\bullet$ be two minimal perfect complexes of free $A$-modules that are quasiisomorphic, i.e.\ there is a  quasiisomorphism $P^\bullet\to Q^\bullet$. Tensoring with $A/\fm$ yields that
$$
P^i/\fm P^i\cong H^i(P^\bullet\otimes_AA/\fm)\cong H^i(Q^\bullet\otimes_AA/\fm )\cong Q^i/\fm Q^i 
$$ 
for all $i$. Thus $P^i\to Q^i$ is an isomorphism as required. 
\eproof

In the following let us call a perfect complex {\em non-trivial} if it is not quasiisomorphic to the zero complex.

\bproof[Proof of Proposition \ref{1.5}]
Suppose that the complex $\F_n^\bullet$ as in \eqref{eq1} is quasiisomorphic to a direct sum of  two perfect free complexes $P^\bullet$ and $Q^\bullet$. 
By Lemma \ref{1.6}(1) we may suppose that both  $P^\bullet$ and $Q^\bullet$ are minimal. According to Lemma \ref{1.6}(2) we get a splitting 
$$
\F_n^\bullet\cong P^\bullet\oplus Q^\bullet. 
$$
Now $\F_n^0\cong A$ is an indecomposable $A$-module, hence either $P^0\cong A$ or $Q^0\cong A$. We may suppose that $P^0\cong A$. Let us show by  induction on $k$ that that then $Q^{-k}=0$ so that in total $Q^\bullet=0$
and $\F^\bullet_n=P^\bullet$. 

Suppose that $\F_n^i=P^i$ for $i=0,\, -1, \ldots, -k+1$ so that $Q^i=0$ in the same range. Then necessarily the map $\partial:Q^{-k}\to \F_n^{-k+1}$ is the zero map. However, inspecting the matrix $M$, the map 
$$
\bar\partial :\F_n^{-k}/(x,y) \F_n^{-k}\to (x,y)^2 \F_n^{-k+1}/(x,y)^3 \F_n^{-k+1}
$$
induced by the differential is injective. As $\partial$ vanishes on the direct summand $Q^{-k}$ of $\F_n^{-k}$, this forces $Q^{-k}=0$. 
\eproof

\bcor\label{1.7}
Theorem \ref{1.3} is true for local rings, i.e.\ a local noetherian ring $A$ is strongly regular if and only if it is a 1-dimensional regular local ring or a field. 
\ecor

\bproof
The `if'-part is clear from Remark \ref{1.2}(1)  above. To show the converse, note first that a strongly regular local ring $A$ is regular by Remark \ref{1.2}(3). Suppose that $\dim A\ge 2$. In that case $A$ being regular, $\fm_A$ contains a regular sequence $(x,y)$. By Proposition \ref{1.5} there exists a perfect complex  $\F^\bullet_n$ of length $n=r+1$ which is not quasiisomorphic to a direct sum of perfect complexes of length $<n$.   Hence $A$ cannot be $r$-regular for any $r\ge 0$.
\eproof

To deduce Theorem \ref{1.3} for arbitrary noetherian rings we need the following observation.

\blem\label{1.8}
Strongly $r$-regular is a local property. More precisely, the following are equivalent. 
\bnum[(1)]
\item $A$ is strongly $r$-regular. 

\item Every localization $A_\fp$ with $\fp\in \Spec A$ is strongly $r$-regular. 

\item Every localization $A_\fm$ with respect to a maximal ideal is strongly $r$-regular. 
\enum
The same is true for the property of being strongly regular. 
\elem 

\bproof
We first treat the property of strong $r$-regularity. 
(2)$\Rightarrow$(3) is trivial. To show (1)$\Rightarrow$(2) let $\fp\subseteq A$ be a prime of $A$ and let $P^\bullet$ be a perfect complex over $A_\fp$. Multiplying its differentials with suitable elements in $A\backslash \fp$ it is easily seen that there is a perfect complex $Q^\bullet $ over $A$ such that $P^\bullet\cong (Q^{\bullet})_\fp$. Since by assumption $Q^\bullet$ is a direct sum of perfect complexes of length $\le r$ the same is true for $P^\bullet$. 

To show (3)$\Rightarrow$(1), let $A$ be a ring for which $A_\fm$ is strongly $r$-regular for every maximal ideal $\fm$ of $A$. By Corollary \ref{1.7} $A_\fm$ is a field or is a discrete valuation ring. Hence $A$ can be written as a product $A\cong \prod_{i=1}^s A_i$, where $A_i$ is a Dedekind domain or a field, and so $A$ is $1$-regular by Remark  \ref{1.2}(1). 

Let us finally consider the property strongly regular. If $A$ is strongly regular then it is strongly $r$-regular for some $r$ and so are all localizations $A_\fp$, $\fp\in \Spec A$, as shown above. The implication (3)$\Rightarrow$(1) for the property `strongly regular' follows finally with the same arguments as before. 
\eproof

\bproof[Proof of Theorem \ref{1.3}]
The `if'-part follows from Remark \ref{1.2}. To prove the converse let $A$ be a strongly regular ring. Every noetherian ring has a decomposition $A= \prod_{i=1}^s A_i$ as a direct product of rings, where $\Spec A_i$, $i=1,\ldots s,$ are just the connected components of $\Spec A$, see e.g.\  \cite[Chapt.\ 1, Ex.\ 22]{AM}. Thus, we may suppose that $\Spec A$ is connected. According to Lemma \ref{1.8} for every maximal ideal $\fm\subseteq A$ the localization  $A_\fm$ is again strongly regular. 
Applying Corollary \ref{1.5}, $A_\fm$ is either a field or a discrete valuation domain. Hence $A$  is either a field or a Dedekind domain.
\eproof

It is interesting to note that with the reasoning in Proposition \ref{1.5} one can also construct many further indecomposable complexes. The construction is as follows.

\bsit\label{1.9}Let $(A,\fm)$ be a local noetherian ring. 
Let us call a complex of free $A$-modules
\bdi
F^\bullet:\quad 0 & \rTo & F^{-d} &\rTo^{\partial} &F^{1-d}  & \rTo^{\partial} & \cdots &\rTo^{\partial} & F^{-2} & \rTo^{\partial} & F^{-1} &\rTo^{\partial} &F^0&\rTo &  0
\edi
{\em irreducible} if 
\bnum[(a)]
\item $F^0\cong A$, and 

\item  for every $i\ge 0$ there is an ideal $\fa\subseteq \fm$ and $s\ge 1$ such that 
$\partial(F^i)\subseteq \fa^sF^i$ and the induced map
$$
\bar \partial: F^{i-1}/\fa F^{i-1}\to \fa^sF^i/\fa^{s+1} F^i
$$ 
is injective. 
\enum
\esit

Using the same reasoning as in the proof of Proposition \ref{1.5} we have the following result. 

\blem\label{1.10}
If  $F^\bullet $ is irreducible then $F^\bullet$ is not isomorphic to a direct sum of two non-trivial perfect complexes. In particular $F^\bullet$ is indecomposable.\hfill $\Box$
\elem

\bexa
(1) Let $(A,\fm)$ be a local ring and let $x_1,\ldots, x_n$ be a regular sequence in $A$. We consider the Koszul complex $K^\bullet:=K^\bullet(x_1,\ldots, x_n)$. More explicitly
$$
K^\bullet: \quad 0\to K^{-n}\to K^{-n+1}\to\ldots\to K^{-1}\to K^0\to 0\,, 
$$
where $K^{-p}=\bigwedge^p A^n$ is freely generated by all wedge products $e_{i_1}\wedge \ldots \wedge e_{i_p}$ with $1\le i_1<\ldots < i_p\le n$, where $e_1,\ldots,e_n$ denotes the standard basis of $A^n$. 
It is easy to see that $K^\bullet$ satisfies conditions (a), (b) in \ref{1.9} with $\fa=(x_1,\ldots, x_n)$ and $s=1$. Using Lemma \ref{1.10} the Koszul complex is indecomposable. 

(2) Let $x_1,\ldots, x_n$ and  $K^\bullet$ as in (1). Composing the two natural maps $K^\bullet[-n]  \to K^{-n}$ and $K^{-n}\cong K^0\hto  K^\bullet$ yields a morphism of complexes $h: K^\bullet[-n] \to K^\bullet$. Its cone $C^\bullet=C^\bullet(h)$  sits in an exact sequence of complexes
$$
0\to K^\bullet\to C^\bullet\to K^\bullet[-n+1]\to 0\,.
$$
In particular $C^k=K^{k-n+1}\oplus K^k$ for all $k$. 
Dividing out the acyclic subcomplex $0\to K^{-n}\oplus 0\to \partial_{C^\bullet}(K^{-n}\oplus 0)\to 0$ of $C^\bullet(h)$ yields a complex $F^\bullet$. The differential 
$$
\partial: F^{-1}\cong K^{-1}\to F^0\cong K^{-n+1}
$$
is up to a sign explicitly given by 
$$
e_i\mapsto \sum_{j=1}^n x_ix_j \cdot e_1\wedge \ldots\wedge\hat e_j\wedge \ldots\wedge e_n. 
$$
In the case $n=2$ the differential is basically given by the matrix $M$ in \ref{1.4}. We call $F^\bullet$ the {\em iterated Koszul complex}. 

Using Lemma \ref{1.10} the iterated Koszul complex is not isomorphic to a direct sum of two non-trivial perfect complexes. In particular it is indecomposable. 

We note that one can repeat the construction to the left or to the right to obtain multiply iterated Koszul complexes, and these complexes are  arbitrarily long and again indecomposable.  
\eexa

\section{Noetherian schemes}

The schemes $X$ that we consider are all assumed to be separated. 
 We recall that a complex
$\cC^\bullet$ of $\cO_X$-modules is  {\em perfect} (cf.\  \cite{Il}) if 
\bnum[(a)]
\item $\cC^\bullet$ has bounded coherent cohomology, and  

\item for every point $P\in X$ the stalk $\cC^\bullet_P$ is quasiisomorphic to a finite complex $F^\bullet$ of  finite free $\cO_{X,P}$-modules. 
\enum
The minimal length of such a complex $F^\bullet$ is called the {\em width} of $\cC^\bullet$ at $P$ and denoted by $\wid_P(\cC^\bullet)$. Moreover we set $\wid(\cC^\bullet):=\max_{P\in X} \wid_P (\cC^\bullet).$ 

We say that $X$ has {\em strong  global dimension $\le r$} if every perfect complex $\cC^\bullet$ of $\cO_X$-modules is quasiisomorphic to a direct sum $\bigoplus_{i=1}^s \cC_{(i)}^\bullet$ of perfect complexes of width $\le r$. 

\brem\label{2.1}
1) Every perfect complex on a scheme is quasiisomorphic to a finite complex of quasicoherent $\cO_X$-modules. 

2) To the best of our knowledge it is not known whether every perfect complex on a noetherian scheme $X$ is quasiisomorphic to a bounded complex of locally free coherent modules.  A closely related question is whether every coherent sheaf on $X$ is isomorphic to a quotient of a locally free coherent sheaf. If $X$ admits an ample line bundle or if it is locally factorial then this is indeed true. For the problem in general see \cite{Gr, Sch}.
\erem

\bprop\label{2.2}
If $X$ is noetherian and has strong global dimension $\le r$ then $X$ is a regular scheme of dimension $\le r$. 
\eprop 

In order to deduce this proposition we need the following result, see Lemma 2.6 in \cite{Ne}.

\bthm\label{2.3}
Let X be a quasicompact, separated scheme and let $U \subseteq X$ be a quasicompact, open subscheme. Let $\cF^\bullet$ be a complex of quasicoherent
$\cO_X$-modules, and let $\cP^\bullet$ be a perfect complex on $U$. Suppose that we are given a map in the derived category of $U$ of the form $p:\cP^\bullet\to\cF^\bullet|U$. Then there exists a perfect complex $\cQ^\bullet$ on $U$  and a morphism of complexes $q:\cQ^\bullet\to \cF^\bullet$ so that the map
\be\label{eq2.2}
p\oplus q: \cP^\bullet\oplus \cQ^\bullet\to \cF^\bullet|U
\ee
lifts to a morphism on $X$, i.e.\ there exists a perfect complex $\cG^\bullet$ on $X$ restricting to $\cP^\bullet\oplus \cQ^\bullet$ on $U$, and  a morphism  $\alpha: \cG^\bullet\to\cF^\bullet$ defined on X, which restricts on U to the given map \eqref{eq2.2}.
\ethm

\bproof[Proof of Proposition \ref{2.2}] 
Let $P\in X$ be a closed point and $U=\Spec A \subseteq X$ an open affine neighborhood.
On $U$ we chose a free resolution $\cF^\bullet$ of the residue field $k_P=\cO_X/\fm_P$, where $\fm_P$ denotes the maximal ideal sheaf of $P$. Shrinking $U$ we may assume that $\cF^\bullet$ is the minimal resolution at $P$ in degrees $\ge -r-1$. Clearly $\cF^\bullet$ is quasiisomorphic to $k_P$ and so can be considered as an element of the derived category of $X$. Cutting this resolution at place $-r-1$ yields a subcomplex of $\cF^\bullet$ on $U$
$$
\cP^\bullet:\quad 0\to \cP^{-r-1}\to \ldots \to\cP^{-1}\to\cP^0\to 0\,
$$
so that $\cP^i=\cF^i$ for $i\ge-r-1$ and $\cP^i=0$ otherwise. According to Theorem \ref{2.3} there is a perfect complex $\cQ^\bullet$ on $U$ and a morphism $q:\cQ^\bullet \to \cF^\bullet$ on $U$ such that the map
\be\label{eq2.3}
p\oplus q:  \cP^\bullet\oplus \cQ^\bullet\to \cF^\bullet\cong k_P|U
\ee
lifts to a morphism $\alpha: \cG^\bullet\to k_P$ on $X$, where $\cG^\bullet$ is perfect on $X$; here $p:\cP^\bullet\to \cF^\bullet$ denotes the inclusion map. 

Assume now that $X$ is $r$-regular. Then $\cG^\bullet$ can be decomposed into a direct sum  $\bigoplus_{i=1}^s \cG^\bullet_{(i)}$, where $\cG^\bullet_{(i)}$ is a perfect complex of width $\le r$ on $X$. Restricting $\alpha$ to direct summands gives the natural map $p:\cP^\bullet\to k_P$,  and another map $\alpha_i: \cG^\bullet_{(i)}\to k_P$. 
There is an index $i$ such that the map $\alpha_i:H^0(\cG^\bullet_{(i)})\to k_P$ is onto. 

Projecting $\cG^\bullet$ onto $\cG^\bullet_{(i)}$ yields a map $\pi_i: \cP^\bullet\to\cG^\bullet_{(i)}|U$, and similarly projecting $\cG^\bullet|U$ onto $\cP^\bullet$ yields a morphisms  $p_i: \cG^\bullet_{(i)}|U\to \cP^\bullet$,
which is the restriction of $p$.

In the derived category of $U$ the diagrams
$$\bdi
\cG^\bullet_{(i)}|U&&\rTo^{p_i} && \cP^\bullet\\
&\rdTo<{\alpha_i} && \ldTo>{p}\\
&& k_P
\edi
\quad\mbox{and}\quad 
\bdi
\cP^\bullet&&\rTo^{\pi_i} && \cG^\bullet_{(i)}|U\\
&\rdTo<{p}  && \ldTo>{\alpha_i}\\
&& k_P
\edi
$$
commute up to multiplication by a non-zero element of $k_P$. 
Thus also the diagram on $U$
$$
\bdi
\cP^\bullet&&\rTo^{p_i\circ \pi_i} && \cP^\bullet\\
&\rdTo<{p}  && \ldTo>{p}\\
&& k_P
\edi
$$
commutes up to multiplication by a non-zero element of $k_P$. Since $\cG^\bullet_{(i)}|U$ has length $\le r$ it follows that 
$p_i\circ \pi_i$ is the zero map on 
$\cP^{-r-1}$.
On the other hand, since $\cF^\bullet_P$ is the minimal resolution of $k_P$ at $P$ in degrees $\ge -r-1$ the map $p_i\circ \pi_i$ induces an isomorphism of complexes, which is only possible if $\cP^{-r-1}=0$. Hence $\pd(k_P)\le r$ as desired.
\eproof

\bexa\label{2.4}
A regular scheme $X$ of dimension $\le 1$ is of strong global dimension $\le 1$. 

\bproof
The proof is similar to that of Remark \ref{1.2}(1).
Again it suffices to show that every perfect complex $\cF^\bullet$ on $X$ is formal. We may assume that $\cF^i$ vanishes for $i> 0$ and consider its cohomology sheaf $\cH:=H^0(\cF^\bullet)$ as a complex concentrated in degree 0. 
Thus, we obtain an exact sequence of complexes 
\be\label{eq2.1}
0\to \cG^\bullet:=\ker(\cF^\bullet\to\cH)\to\cF^\bullet\to \cH\to 0\,.
\ee
By the same spectral sequence argument as in \ref{1.2}(1) the group
$\Ext^{1}_X(\cH, \cG^\bullet)$ vanishes. Hence the sequence \eqref{eq2.1} splits, i.e.\ $\cF^\bullet$ is quasiisomorphic to $ \cG^\bullet\oplus\cH$ . A simple induction  now yields that $\cF^\bullet$ is formal, as required.
\eproof
\eexa

We have the following result. 

\bthm\label{2.5}
A noetherian scheme $X$ is of finite strong global dimension if and only if it is regular of dimension $\le 1$. 
\ethm

\bproof
By Proposition \ref{2.2} we know that $X$ is a regular scheme. 
Assume that $X$ is  of dimension $\ge 2$. We may assume that $X$ is connected. 
Let $D$ and $E$ be two effective divisors on $X$ which intersect in a non-empty subset of codimension $\ge 2$. We consider the $\cO_X$-modules 
$$
\cF^0:=\cO_X \quad \mbox{ and } \quad \cF^n:=\cO_X(-D+n(D+E))\oplus \cO_X(-E+n(D+E)) \quad \mbox{for }n\le -1\,.
$$
For any $\cO_X$-module $\cM$ there are natural maps 
$$
x: \cM(-D)\to \cM\quad\mbox{and}\quad
y: \cM(-E)\to \cM\,.
$$
Combined into the matrix
$$
M:=\bmat xy &-x^2 \\ y^2& -xy
\emat 
$$
this yields a map $\partial^i: \cF^{i-1}\to \cF^{i}$, $i\le - 1$, with  
$\partial^{i}\circ \partial^{i-1}=0$. Moreover  let $\partial^0:\cF^{-1}\to \cF^0$ be the map given by the matrix $\bsmat x\\-y\esmat$. With these maps we obtain for any $n\ge 0$ a perfect complex  
\bdi
\cF^\bullet:\quad 0& \rTo &\cF^{-n} &\rTo^{\partial^{-n+1}} &
\cF^{-n+1}&\rTo^{\partial^{-n+2}} &\ldots &\rTo^{\partial^{-1}} &
\cF^{-1}& \rTo^{\partial^{0}} & \cF^0 &\rTo & 0\,.
\edi 
Assume  $\cF^\bullet\cong \bigoplus_i \cF^\bullet_{(i)}$ were a direct sum of perfect complexes of length $<n$. 
At a point $P\in D\cap E$ in the intersection, the stalk $\cF^\bullet_P$ were then a direct sum of perfect complexes over $\cO_{X,P}$ of length $<n$, contradicting Proposition \ref{1.5}.
\eproof

\end{document}